\documentclass{amsart}
\usepackage{amsfonts}

\usepackage{graphicx}
\usepackage{amscd}

\newtheorem{theorem}{Theorem}
\theoremstyle{plain}

\newtheorem{corollary}{Corollary}

\newtheorem{remark}{Remark}

\numberwithin{equation}{section}

\input{tcilatex}

\begin{document}
\title[Ostrowski Inequality]{An Ostrowski Type Inequality for Convex Functions}
\author{S.S. Dragomir}
\address{School of Communications and Informatics\\
Victoria University of Technology\\
PO Box 14428\\
Melbourne City MC\\
Victoria 8001, Australia.}
\email{sever@matilda.vu.edu.au}
\urladdr{http://rgmia.vu.edu.au/SSDragomirWeb.html}
\date{June 13, 2001.}
\subjclass{Primary 26D14, 26D99.}
\keywords{Ostrowski Inequality, Hermite-Hadamard inequality, Integral Means,
Probability density function, Divergence measures.}

\begin{abstract}
An Ostrowski type integral inequality for convex functions and applications
for quadrature rules and integral means are given. A refinement and a
counterpart result for Hermite-Hadamard inequalities are obtained and some
inequalities for pdf's and $\left( HH\right) -$divergence measure are also
mentioned.
\end{abstract}

\maketitle

\section{Introduction}

The following result is known in the literature as Ostrowski's inequality 
\cite{1b}.

\begin{theorem}
\label{ta}Let $f:\left[ a,b\right] \rightarrow \mathbb{R}$ be a
differentiable mapping on $\left( a,b\right) $ with the property that $%
\left| f^{\prime }\left( t\right) \right| \leq M$ for all $t\in \left(
a,b\right) $. Then 
\begin{equation}
\left| f\left( x\right) -\frac{1}{b-a}\int_{a}^{b}f\left( t\right) dt\right|
\leq \left[ \frac{1}{4}+\frac{\left( x-\frac{a+b}{2}\right) ^{2}}{\left(
b-a\right) ^{2}}\right] \left( b-a\right) M  \label{1.1}
\end{equation}
for all $x\in \left[ a,b\right] $.\newline
The constant $\frac{1}{4}$ is the best possible in the sense that it cannot
be replaced by a smaller constant.
\end{theorem}

A simple proof of this fact can be done by using the identity: 
\begin{equation}
f\left( x\right) =\frac{1}{b-a}\int_{a}^{b}f\left( t\right) dt+\frac{1}{b-a}%
\int_{a}^{b}p\left( x,t\right) f^{\prime }\left( t\right) dt,\;\;x\in \left[
a,b\right] ,  \label{1.2}
\end{equation}
where 
\begin{equation*}
p\left( x,t\right) :=\left\{ 
\begin{array}{lll}
t-a & \text{if} & a\leq t\leq x \\ 
&  &  \\ 
t-b & \text{if} & x<t\leq b
\end{array}
\right.
\end{equation*}
which holds for absolutely continuous functions $f:\left[ a,b\right]
\rightarrow \mathbb{R}$.

The following Ostrowski type result holds (see \cite{2b}, \cite{3b} and \cite
{4b}).

\begin{theorem}
\label{tb}Let $f:\left[ a,b\right] \rightarrow \mathbb{R}$ be absolutely
continuous on $\left[ a,b\right] $. Then, for all $x\in \left[ a,b\right] $,
we have: 
\begin{eqnarray}
&&\left| f\left( x\right) -\frac{1}{b-a}\int_{a}^{b}f\left( t\right)
dt\right|  \label{1.3} \\
&\leq &\left\{ 
\begin{array}{lll}
\left[ \frac{1}{4}+\left( \frac{x-\frac{a+b}{2}}{b-a}\right) ^{2}\right]
\left( b-a\right) \left\| f^{\prime }\right\| _{\infty } & \text{if} & 
f^{\prime }\in L_{\infty }\left[ a,b\right] ; \\ 
&  &  \\ 
\frac{1}{\left( p+1\right) ^{\frac{1}{p}}}\left[ \left( \frac{x-a}{b-a}%
\right) ^{p+1}+\left( \frac{b-x}{b-a}\right) ^{p+1}\right] ^{\frac{1}{p}%
}\left( b-a\right) ^{\frac{1}{p}}\left\| f^{\prime }\right\| _{q} & \text{if}
& f^{\prime }\in L_{q}\left[ a,b\right] , \\ 
&  & \frac{1}{p}+\frac{1}{q}=1,\;p>1; \\ 
\left[ \frac{1}{2}+\left| \frac{x-\frac{a+b}{2}}{b-a}\right| \right] \left\|
f^{\prime }\right\| _{1}; &  & 
\end{array}
\right.  \notag
\end{eqnarray}
where $\left\| \cdot \right\| _{r}$ \ ($r\in \left[ 1,\infty \right] $) are
the usual Lebesgue norms on $L_{r}\left[ a,b\right] $, i.e., 
\begin{equation*}
\left\| g\right\| _{\infty }:=ess\sup\limits_{t\in \left[ a,b\right] }\left|
g\left( t\right) \right|
\end{equation*}
and 
\begin{equation*}
\left\| g\right\| _{r}:=\left( \int_{a}^{b}\left| g\left( t\right) \right|
^{r}dt\right) ^{\frac{1}{r}},\;r\in \lbrack 1,\infty ).
\end{equation*}
The constants $\frac{1}{4}$, $\frac{1}{\left( p+1\right) ^{\frac{1}{p}}}$
and $\frac{1}{2}$ respectively are sharp in the sense presented in Theorem 
\ref{ta}.
\end{theorem}

The above inequalities can also be obtained from Fink's result in \cite{5b}
on choosing $n=1$ and performing some appropriate computations.

If one drops the condition of absolute continuity and assumes that $f$ is
H\"{o}lder continuous, then one may state the result (see \cite{6b}):

\begin{theorem}
\label{tc}Let $f:\left[ a,b\right] \rightarrow \mathbb{R}$ be of $r-H-$%
H\"{o}lder type, i.e., 
\begin{equation}
\left| f\left( x\right) -f\left( y\right) \right| \leq H\left| x-y\right|
^{r},\;\text{for all \ }x,y\in \left[ a,b\right] ,  \label{1.4}
\end{equation}
where $r\in (0,1]$ and $H>0$ are fixed. Then for all $x\in \left[ a,b\right] 
$ we have the inequality: 
\begin{eqnarray}
&&\left| f\left( x\right) -\frac{1}{b-a}\int_{a}^{b}f\left( t\right)
dt\right|  \label{1.5} \\
&\leq &\frac{H}{r+1}\left[ \left( \frac{b-x}{b-a}\right) ^{r+1}+\left( \frac{%
x-a}{b-a}\right) ^{r+1}\right] \left( b-a\right) ^{r}.  \notag
\end{eqnarray}
The constant $\frac{1}{r+1}$ is also sharp in the above sense.
\end{theorem}

Note that if $r=1$, i.e., $f$ is Lipschitz continuous, then we get the
following version of Ostrowski's inequality for Lipschitzian functions (with 
$L$ instead of $H$) (see \cite{7b}) 
\begin{equation}
\left| f\left( x\right) -\frac{1}{b-a}\int_{a}^{b}f\left( t\right) dt\right|
\leq \left[ \frac{1}{4}+\left( \frac{x-\frac{a+b}{2}}{b-a}\right) ^{2}\right]
\left( b-a\right) L.  \label{1.6}
\end{equation}
Here the constant $\frac{1}{4}$ is also best.

Moreover, if one drops the continuity condition of the function, and assumes
that it is of bounded variation, then the following result may be stated
(see \cite{8b}).

\begin{theorem}
\label{td}Assume that $f:\left[ a,b\right] \rightarrow \mathbb{R}$ is of
bounded variation and denote by $\bigvee\limits_{a}^{b}\left( f\right) $ its
total variation. Then 
\begin{equation}
\left| f\left( x\right) -\frac{1}{b-a}\int_{a}^{b}f\left( t\right) dt\right|
\leq \left[ \frac{1}{2}+\left| \frac{x-\frac{a+b}{2}}{b-a}\right| \right]
\bigvee\limits_{a}^{b}\left( f\right)  \label{1.7}
\end{equation}
for all $x\in \left[ a,b\right] $.\newline
The constant $\frac{1}{2}$ is the best possible.
\end{theorem}

If we assume more about $f$, i.e., $f$ is monotonically increasing, then the
inequality (\ref{1.7}) may be improved in the following manner \cite{9b}
(see also \cite{10b}).

\begin{theorem}
\label{te}Let $f:\left[ a,b\right] \rightarrow \mathbb{R}$ be monotonic
nondecreasing. Then for all $x\in \left[ a,b\right] $, we have the
inequality: 
\begin{eqnarray}
&&\left| f\left( x\right) -\frac{1}{b-a}\int_{a}^{b}f\left( t\right)
dt\right|  \label{1.8} \\
&\leq &\frac{1}{b-a}\left\{ \left[ 2x-\left( a+b\right) \right] f\left(
x\right) +\int_{a}^{b}sgn\left( t-x\right) f\left( t\right) dt\right\} 
\notag \\
&\leq &\frac{1}{b-a}\left\{ \left( x-a\right) \left[ f\left( x\right)
-f\left( a\right) \right] +\left( b-x\right) \left[ f\left( b\right)
-f\left( x\right) \right] \right\}  \notag \\
&\leq &\left[ \frac{1}{2}+\left| \frac{x-\frac{a+b}{2}}{b-a}\right| \right] %
\left[ f\left( b\right) -f\left( a\right) \right] .  \notag
\end{eqnarray}
All the inequalities in (\ref{1.8}) are sharp and the constant $\frac{1}{2}$
is the best possible.
\end{theorem}

In this paper we establish an Ostrowski type inequality for convex
functions. Applications for quadrature rules, for integral means, for
probability distribution functions, and for $HH-$divergences in Information
Theory are also considered.

\section{The Results}

The following theorem providing a lower bound for the Ostrowski difference $%
\int_{a}^{b}f\left( t\right) dt-\left( b-a\right) f\left( x\right) $ holds.

\begin{theorem}
\label{t1}Let $f:\left[ a,b\right] \rightarrow \mathbb{R}$ be a convex
function on $\left[ a,b\right] $. Then for any $x\in \left( a,b\right) $ we
have the inequality: 
\begin{equation}
\frac{1}{2}\left[ \left( b-x\right) ^{2}f_{+}^{\prime }\left( x\right)
-\left( x-a\right) ^{2}f_{-}^{\prime }\left( x\right) \right] \leq
\int_{a}^{b}f\left( t\right) dt-\left( b-a\right) f\left( x\right) .
\label{2.1}
\end{equation}
The constant $\frac{1}{2}$\ in the left hand side of (\ref{2.1}) is sharp in
the sense that it cannot be replaced by a larger constant.
\end{theorem}

\begin{proof}
It is easy to see that for any locally absolutely continuous function $%
f:\left( a,b\right) \rightarrow \mathbb{R}$, we have the identity 
\begin{equation}
\int_{a}^{x}\left( t-a\right) f^{\prime }\left( t\right)
dt+\int_{x}^{b}\left( t-b\right) f^{\prime }\left( t\right) dt=f\left(
x\right) -\int_{a}^{b}f\left( t\right) dt,  \label{2.2}
\end{equation}
for any $x\in \left( a,b\right) $ where $f^{\prime }$ is the derivative of $%
f $ which exists a.e. on $\left( a,b\right) .$

Since $f$ is convex, then it is locally Lipschitzian and thus (\ref{2.2})
holds. Moreover, for any $x\in \left( a,b\right) ,$ we have the inequalities 
\begin{equation}
f^{\prime }\left( t\right) \leq f_{-}^{\prime }\left( x\right) \text{ for
a.e. }t\in \left[ a,x\right]  \label{2.3}
\end{equation}
and 
\begin{equation}
f^{\prime }\left( t\right) \geq f_{+}^{\prime }\left( x\right) \text{ for
a.e. }t\in \left[ x,b\right] .  \label{2.4}
\end{equation}
If we multiply (\ref{2.3}) by $t-a\geq 0,$ $t\in \left[ a,x\right] ,$ and
integrate on $\left[ a,x\right] $, we get 
\begin{equation}
\int_{a}^{x}\left( t-a\right) f^{\prime }\left( t\right) dt\leq \frac{1}{2}%
\left( x-a\right) ^{2}f_{-}^{\prime }\left( x\right)  \label{2.5}
\end{equation}
and if we multiply (\ref{2.4}) by $b-t\geq 0,$ $t\in \left[ x,b\right] ,$
and integrate on $\left[ x,b\right] ,$ we also have 
\begin{equation}
\int_{x}^{b}\left( b-t\right) f^{\prime }\left( t\right) dt\geq \frac{1}{2}%
\left( b-x\right) ^{2}f_{+}^{\prime }\left( x\right) .  \label{2.6}
\end{equation}
Finally, if we subtract (\ref{2.6}) from (\ref{2.5}) and use the
representation (\ref{2.2}) we deduce the desired inequality (\ref{2.1}).

Now, assume that (\ref{2.1}) holds with a constant $C>0$ instead of $\frac{1%
}{2},$ i.e., 
\begin{equation}
C\left[ \left( b-x\right) ^{2}f_{+}^{\prime }\left( x\right) -\left(
x-a\right) ^{2}f_{-}^{\prime }\left( x\right) \right] \leq
\int_{a}^{b}f\left( t\right) dt-\left( b-a\right) f\left( x\right) .
\label{2.7}
\end{equation}
Consider the convex function $f_{0}\left( t\right) :=k\left| t-\frac{a+b}{2}%
\right| $, $k>0,$ $t\in \left[ a,b\right] .$ Then 
\begin{equation*}
f_{0^{+}}^{\prime }\left( \frac{a+b}{2}\right) =k,\;\;\;f_{0^{-}}^{\prime
}\left( \frac{a+b}{2}\right) =-k,\;\;\;f_{0}\left( \frac{a+b}{2}\right) =0
\end{equation*}
and 
\begin{equation*}
\int_{a}^{b}f_{0}\left( t\right) dt=\frac{1}{4}k\left( b-a\right) ^{2}.
\end{equation*}
If in (\ref{2.7}) we choose $f_{0}$ as above and $x=\frac{a+b}{2},$ then we
get 
\begin{equation*}
C\left[ \frac{1}{4}\left( b-a\right) ^{2}k+\frac{1}{4}\left( b-a\right) ^{2}k%
\right] \leq \frac{1}{4}k\left( b-a\right) ^{2},
\end{equation*}
which gives $C\leq \frac{1}{2},$ and the sharpness of the constant is proved.
\end{proof}

Now, recall that the following inequality, which is well known in the
literature as the \textit{Hermite-Hadamard inequality} for convex functions,
holds: 
\begin{equation}
f\left( \frac{a+b}{2}\right) \leq \frac{1}{b-a}\int_{a}^{b}f\left( t\right)
dt\leq \frac{f\left( a\right) +f\left( b\right) }{2}.  \tag{HH}  \label{HH}
\end{equation}

The following corollary which improves the first Hermite-Hadamard inequality
(\ref{HH}) holds.

\begin{corollary}
\label{c1}Let $f:\left[ a,b\right] \rightarrow \mathbb{R}$ be a convex
function on $\left[ a,b\right] $. Then 
\begin{eqnarray}
0 &\leq &\frac{1}{8}\left[ f_{+}^{\prime }\left( \frac{a+b}{2}\right)
-f_{-}^{\prime }\left( \frac{a+b}{2}\right) \right] \left( b-a\right)
\label{2.8} \\
&\leq &\frac{1}{b-a}\int_{a}^{b}f\left( t\right) dt-f\left( \frac{a+b}{2}%
\right) .  \notag
\end{eqnarray}
The constant $\frac{1}{8}$ is sharp.
\end{corollary}

The proof is obvious by the above theorem. The sharpness of the constant is
obtained for $f_{0}\left( t\right) :=k\left| t-\frac{a+b}{2}\right| ,$ $t\in %
\left[ a,b\right] ,$ $k>0.$

When $x$ is a point of differentiability, we may state the following
corollary as well.

\begin{corollary}
\label{c2}Let $f$ be as in Theorem \ref{t1}. If $x\in \left( a,b\right) $ is
a point of differentiability for $f$, then 
\begin{equation}
\left( \frac{a+b}{2}-x\right) f^{\prime }\left( x\right) \leq \frac{1}{b-a}%
\int_{a}^{b}f\left( t\right) dt-f\left( x\right) .  \label{2.9}
\end{equation}
\end{corollary}

\begin{remark}
\label{r1}If $f:I\subseteq \mathbb{R\rightarrow R}$ is convex on $I$ and if
we choose $x\in $\r{I} \ (\r{I} is the interior of $I$), $b=x+\frac{h}{2}$, $%
a=x-\frac{h}{2}$, $h>0$ is such that $a,b\in I$, then from (\ref{2.1}) we
may write 
\begin{equation}
0\leq \frac{1}{8}h^{2}\left[ f_{+}^{\prime }\left( x\right) -f_{-}^{\prime
}\left( x\right) \right] \leq \int_{x-\frac{h}{2}}^{x+\frac{h}{2}}f\left(
t\right) dt-hf\left( x\right) ,  \label{2.10}
\end{equation}
and the constant $\frac{1}{8}$ is sharp in (\ref{2.10}).
\end{remark}

The following result providing an upper bound for the Ostrowski difference $%
\int_{a}^{b}f\left( t\right) dt-\left( b-a\right) f\left( x\right) $ also
holds.

\begin{theorem}
\label{t2}Let $f:\left[ a,b\right] \rightarrow \mathbb{R}$ be a convex
function on $\left[ a,b\right] .$ Then for any $x\in \left[ a,b\right] ,$ we
have the inequality: 
\begin{equation}
\int_{a}^{b}f\left( t\right) dt-\left( b-a\right) f\left( x\right) \leq 
\frac{1}{2}\left[ \left( b-x\right) ^{2}f_{-}^{\prime }\left( b\right)
-\left( x-a\right) ^{2}f_{+}^{\prime }\left( a\right) \right] .  \label{2.11}
\end{equation}
The constant $\frac{1}{2}$ is sharp in the sense that it cannot be replaced
by a smaller constant.
\end{theorem}

\begin{proof}
If either $f_{+}^{\prime }\left( a\right) =-\infty $ or $f_{-}^{\prime
}\left( b\right) =+\infty ,$ then the inequality (\ref{2.11}) evidently
holds true.

Assume that $f_{+}^{\prime }\left( a\right) $ and $f_{-}^{\prime }\left(
b\right) $ are finite.

Since $f$ is convex on $\left[ a,b\right] ,$ we have 
\begin{equation}
f^{\prime }\left( t\right) \geq f_{+}^{\prime }\left( a\right) \text{ for
a.e. }t\in \left[ a,x\right]  \label{2.12}
\end{equation}
and 
\begin{equation}
f^{\prime }\left( t\right) \leq f_{-}^{\prime }\left( b\right) \text{ for
a.e. }t\in \left[ x,b\right] .  \label{2.13}
\end{equation}
If we multiply (\ref{2.12}) by $t-a\geq 0,$ $t\in \left[ a,x\right] ,$ and
integrate on $\left[ a,x\right] ,$ then we deduce 
\begin{equation}
\int_{a}^{x}\left( t-a\right) f^{\prime }\left( t\right) dt\geq \frac{1}{2}%
\left( x-a\right) ^{2}f_{+}^{\prime }\left( a\right)  \label{2.14}
\end{equation}
and if we multiply (\ref{2.13}) by $b-t\geq 0,$ $t\in \left[ x,b\right] ,$
and integrate on $\left[ x,b\right] ,$ then we also have 
\begin{equation}
\int_{x}^{b}\left( b-t\right) f^{\prime }\left( t\right) dt\leq \frac{1}{2}%
\left( b-x\right) ^{2}f_{-}^{\prime }\left( b\right) .  \label{2.15}
\end{equation}
Finally, if we subtract (\ref{2.14}) from (\ref{2.15}) and use the
representation (\ref{2.2}), we deduce the desired inequality (\ref{2.11}).

Now, assume that (\ref{2.11}) holds with a constant $D>0$ instead of $\frac{1%
}{2},$ i.e., 
\begin{equation}
\int_{a}^{b}f\left( t\right) dt-\left( b-a\right) f\left( x\right) \leq D 
\left[ \left( b-x\right) ^{2}f_{-}^{\prime }\left( b\right) -\left(
x-a\right) ^{2}f_{+}^{\prime }\left( a\right) \right] .  \label{2.16}
\end{equation}
If we consider the convex function $f_{0}:\left[ a,b\right] \rightarrow 
\mathbb{R}$, $f_{0}\left( t\right) =k\left| t-\frac{a+b}{2}\right| ,$ then
we have $f_{-}^{\prime }\left( b\right) =k$, $f_{+}^{\prime }\left( a\right)
=-k$ and by (\ref{2.16}) we deduce for $x=\frac{a+b}{2}$ that 
\begin{equation*}
\frac{1}{4}k\left( b-a\right) ^{2}\leq D\left[ \frac{1}{4}k\left( b-a\right)
^{2}+\frac{1}{4}k\left( b-a\right) ^{2}\right] ,
\end{equation*}
giving $D\geq \frac{1}{2},$ and the sharpness of the constant is proved.
\end{proof}

The following corollary related to the Hermite-Hadamard inequality is
interesting as well.

\begin{corollary}
\label{c3}Let $f:\left[ a,b\right] \rightarrow \mathbb{R}$ be convex on $%
\left[ a,b\right] $. Then 
\begin{equation}
0\leq \frac{1}{b-a}\int_{a}^{b}f\left( t\right) dt-f\left( \frac{a+b}{2}%
\right) \leq \frac{1}{8}\left[ f_{-}^{\prime }\left( b\right) -f_{+}^{\prime
}\left( a\right) \right] \left( b-a\right)   \label{2.17}
\end{equation}
and the constant $\frac{1}{8}$ is sharp.
\end{corollary}

\begin{remark}
\label{r2}Denote $B:=f_{-}^{\prime }\left( b\right) $, $A:=f_{+}^{\prime
}\left( a\right) $ and assume that $B\neq A,$ i.e., $f$ is not constant on $%
\left( a,b\right) $. Then 
\begin{eqnarray*}
&&\left( b-x\right) ^{2}B-\left( x-a\right) ^{2}A \\
&=&\left( B-A\right) \left[ x-\left( \frac{bB-aA}{B-A}\right) \right] ^{2}-%
\frac{AB}{B-A}\left( b-a\right) ^{2}
\end{eqnarray*}
and by (\ref{2.11}) we get 
\begin{eqnarray}
&&\int_{a}^{b}f\left( t\right) dt-\left( b-a\right) f\left( x\right)
\label{2.18} \\
&\leq &\frac{1}{2}\left( B-A\right) \left\{ \left[ x-\left( \frac{bB-aA}{B-A}%
\right) \right] ^{2}-\frac{AB}{\left( B-A\right) ^{2}}\left( b-a\right)
^{2}\right\}  \notag
\end{eqnarray}
for any $x\in \left[ a,b\right] .$

If $A\geq 0$ then $x_{0}=\frac{bB-aA}{B-A}\in \left[ a,b\right] $ and by (%
\ref{2.18}) we get, choosing $x=\frac{bB-aA}{B-A}$, that 
\begin{equation}
0\leq \frac{1}{2}\frac{AB}{B-A}\left( b-a\right) \leq f\left( \frac{bB-aA}{%
B-A}\right) -\frac{1}{b-a}\int_{a}^{b}f\left( t\right) dt,  \label{2.19}
\end{equation}
which is an interesting inequality in itself.
\end{remark}

\begin{remark}
\label{r3}If $f:I\subseteq \mathbb{R\rightarrow R}$ is convex on $I$ and if
we choose $x\in $\r{I}, $b=x+\frac{h}{2}$, $a=x-\frac{h}{2}$, $h>0$ is such
that $a,b\in I,$ then from (\ref{2.11}) we deduce: 
\begin{equation}
0\leq \int_{x-\frac{h}{2}}^{x+\frac{h}{2}}f\left( t\right) dt-hf\left(
x\right) \leq \frac{1}{8}h^{2}\left[ f_{-}^{\prime }\left( x+\frac{h}{2}%
\right) -f_{+}^{\prime }\left( x-\frac{h}{2}\right) \right] ,  \label{2.20}
\end{equation}
and the constant $\frac{1}{8}$ is sharp.
\end{remark}

\section{The Composite Case}

Consider the division $I_{n}:a=x_{0}<x_{1}<\cdots <x_{n-1}<x_{n}=b$ and
denote $h_{i}:=x_{i+1}-x_{i}$, $i=\overline{0,n-1}.$ If $\xi _{i}\in \left[
x_{i},x_{i+1}\right] $ $\left( i=\overline{0,n-1}\right) $ are intermediate
points, then we will denote by 
\begin{equation}
R_{n}\left( f;I_{n},\mathbf{\xi }\right) :=\sum_{i=0}^{n-1}h_{i}f\left( \xi
_{i}\right)  \label{3.0}
\end{equation}
the Riemann sum associated to $f$, $I_{n}$ and $\mathbf{\xi }$.

The following theorem providing upper and lower bounds for the remainder in
approximating the integral $\int_{a}^{b}f\left( t\right) dt$ of a convex
function $f$ in terms of a general Riemann sum holds.

\begin{theorem}
\label{t3.1}Let $f:\left[ a,b\right] \rightarrow \mathbb{R}$ be a convex
function and $I_{n}$ and $\mathbf{\xi }$ be as above. Then we have: 
\begin{equation}
\int_{a}^{b}f\left( t\right) dt=R_{n}\left( f;I_{n},\mathbf{\xi }\right)
+W_{n}\left( f;I_{n},\mathbf{\xi }\right) ,  \label{3.1}
\end{equation}
where $R_{n}\left( f;I_{n},\mathbf{\xi }\right) $ is the Riemann sum defined
by (\ref{3.0}) and the remainder $W_{n}\left( f;I_{n},\mathbf{\xi }\right) $
satisfies the estimate: 
\begin{eqnarray}
&&\frac{1}{2}\left[ \sum_{i=0}^{n-1}\left( x_{i+1}-\xi _{i}\right)
^{2}f_{+}^{\prime }\left( \xi _{i}\right) -\sum_{i=0}^{n-1}\left( \xi
_{i}-x_{i}\right) ^{2}f_{-}^{\prime }\left( \xi _{i}\right) \right]
\label{3.2} \\
&\leq &W_{n}\left( f;I_{n},\mathbf{\xi }\right)  \notag \\
&\leq &\frac{1}{2}\left[ \left( b-\xi _{n-1}\right) ^{2}f_{-}^{\prime
}\left( b\right) +\sum_{i=1}^{n-1}\right. \left[ \left( x_{i}-\xi
_{i-1}\right) ^{2}f_{-}^{\prime }\left( x_{i}\right) \right.  \notag \\
&&-\left. \left( \xi _{i}-x_{i}\right) ^{2}f_{+}^{\prime }\left(
x_{i}\right) \right] -\left( \xi _{0}-a\right) ^{2}f_{+}^{\prime }\left(
a\right) \bigg].  \notag
\end{eqnarray}
\end{theorem}

\begin{proof}
If we write the inequalities (\ref{2.1}) and (\ref{2.11}) on the interval $%
\left[ x_{i},x_{i+1}\right] $ and for the intermediate points $\xi _{i}\in %
\left[ x_{i},x_{i+1}\right] ,$ then we have 
\begin{eqnarray*}
&&\frac{1}{2}\left[ \left( x_{i+1}-\xi _{i}\right) ^{2}f_{+}^{\prime }\left(
\xi _{i}\right) -\left( \xi _{i}-x_{i}\right) ^{2}f_{-}^{\prime }\left( \xi
_{i}\right) \right]  \\
&\leq &\int_{x_{i}}^{x_{i+1}}f\left( t\right) dt-h_{i}f\left( \xi
_{i}\right)  \\
&\leq &\frac{1}{2}\left[ \left( x_{i+1}-\xi _{i}\right) ^{2}f_{-}^{\prime
}\left( x_{i+1}\right) -\left( \xi _{i}-x_{i}\right) ^{2}f_{+}^{\prime
}\left( x_{i}\right) \right] .
\end{eqnarray*}
Summing the above inequalities over $i$ from $0$ to $n-1$, we deduce 
\begin{eqnarray}
&&\frac{1}{2}\sum_{i=0}^{n-1}\left[ \left( x_{i+1}-\xi _{i}\right)
^{2}f_{+}^{\prime }\left( \xi _{i}\right) -\left( \xi _{i}-x_{i}\right)
^{2}f_{-}^{\prime }\left( \xi _{i}\right) \right]   \label{3.3} \\
&\leq &\int_{a}^{b}f\left( t\right) dt-R_{n}\left( f;I_{n},\mathbf{\xi }%
\right)   \notag \\
&\leq &\frac{1}{2}\left[ \sum_{i=0}^{n-1}\left( x_{i+1}-\xi _{i}\right)
^{2}f_{-}^{\prime }\left( x_{i+1}\right) -\sum_{i=0}^{n-1}\left( \xi
_{i}-x_{i}\right) ^{2}f_{+}^{\prime }\left( x_{i}\right) \right] .  \notag
\end{eqnarray}
However, 
\begin{eqnarray*}
\sum_{i=0}^{n-1}\left( x_{i+1}-\xi _{i}\right) ^{2}f_{-}^{\prime }\left(
x_{i+1}\right)  &=&\left( b-\xi _{n-1}\right) ^{2}f_{-}^{\prime }\left(
b\right) +\sum_{i=0}^{n-2}\left[ \left( x_{i+1}-\xi _{i}\right)
^{2}f_{-}^{\prime }\left( x_{i+1}\right) \right]  \\
&=&\left( b-\xi _{n-1}\right) ^{2}f_{-}^{\prime }\left( b\right)
+\sum_{i=1}^{n-1}\left[ \left( x_{i}-\xi _{i-1}\right) ^{2}f_{-}^{\prime
}\left( x_{i}\right) \right] 
\end{eqnarray*}
and 
\begin{equation*}
\sum_{i=0}^{n-1}\left( \xi _{i}-x_{i}\right) ^{2}f_{+}^{\prime }\left(
x_{i}\right) =\sum_{i=1}^{n-1}\left( \xi _{i}-x_{i}\right) ^{2}f_{+}^{\prime
}\left( x_{i}\right) +\left( \xi _{0}-a\right) ^{2}f_{+}^{\prime }\left(
a\right) 
\end{equation*}
and then, by (\ref{3.3}), we deduce the desired estimate (\ref{3.2}).
\end{proof}

The following corollary may be useful in practical applications.

\begin{corollary}
\label{c3.2}Let $f:\left[ a,b\right] \rightarrow \mathbb{R}$ be a
differentiable convex function on $\left( a,b\right) $. Then we have the
representation (\ref{3.1}) and the remainder $W_{n}\left( f;I_{n},\mathbf{%
\xi }\right) $ satisfies the estimate: 
\begin{eqnarray}
&&\sum_{i=0}^{n-1}\left( \frac{x_{i}+x_{i+1}}{2}-\xi _{i}\right)
h_{i}f^{\prime }\left( \xi _{i}\right)  \label{3.4} \\
&\leq &W_{n}\left( f;I_{n},\xi _{i}\right)  \notag \\
&\leq &\frac{1}{2}\left[ \left( b-\xi _{n-1}\right) ^{2}f_{-}^{\prime
}\left( b\right) -\left( \xi _{0}-a\right) ^{2}f_{+}^{\prime }\left(
a\right) \right]  \notag \\
&&+\sum_{i=1}^{n-1}\left( x_{i}-\frac{\xi _{i}+\xi _{i-1}}{2}\right) \left(
\xi _{i}-\xi _{i-1}\right) f^{\prime }\left( x_{i}\right) .  \notag
\end{eqnarray}
\end{corollary}

We may also consider the mid-point quadrature rule: 
\begin{equation}
M_{n}\left( f,I_{n}\right) :=\sum_{i=0}^{n-1}h_{i}f\left( \frac{x_{i}+x_{i+1}%
}{2}\right) .  \label{3.5}
\end{equation}
Using Corollaries \ref{c1} and \ref{c2}, we may state the following result
as well.

\begin{corollary}
\label{c3.3}Assume that $f:\left[ a,b\right] \rightarrow \mathbb{R}$ is a
convex function on $\left[ a,b\right] $ and $I_{n}$ is a division as above.
Then we have the representation: 
\begin{equation}
\int_{a}^{b}f\left( x\right) dx=M_{n}\left( f,I_{n}\right) +S_{n}\left(
f,I_{n}\right) ,  \label{3.6}
\end{equation}
where $M_{n}\left( f,I_{n}\right) $ is the mid-point quadrature rule given
in (\ref{3.5}) and the remainder $S_{n}\left( f,I_{n}\right) $ satisfies the
estimates: 
\begin{eqnarray}
0 &\leq &\frac{1}{8}\sum_{i=0}^{n-1}\left[ f_{+}^{\prime }\left( \frac{%
x_{i}+x_{i+1}}{2}\right) -f_{-}^{\prime }\left( \frac{x_{i}+x_{i+1}}{2}%
\right) \right] h_{i}^{2}  \label{3.7} \\
&\leq &S_{n}\left( f,I_{n}\right) \leq \frac{1}{8}\sum_{i=0}^{n-1}\left[
f_{-}^{\prime }\left( x_{i+1}\right) -f_{+}^{\prime }\left( x_{i}\right) %
\right] h_{i}^{2}.  \notag
\end{eqnarray}
The constant $\frac{1}{8}$\ is sharp in both inequalities.
\end{corollary}

\section{Inequalities for Integral Means}

We may prove the following result in comparing two integral means.

\begin{theorem}
\label{t4.1}Let $f:\left[ a,b\right] \rightarrow \mathbb{R}$ be a convex
function and $c,d\in \left[ a,b\right] $ with $c<d$. Then we have the
inequalities 
\begin{eqnarray}
&&\frac{a+b}{2}\cdot \frac{f\left( d\right) -f\left( c\right) }{d-c}-\frac{%
df\left( d\right) -cf\left( c\right) }{d-c}+\frac{1}{d-c}\int_{c}^{d}f\left(
x\right) dx  \label{4.1} \\
&\leq &\frac{1}{b-a}\int_{a}^{b}f\left( t\right) dt-\frac{1}{d-c}%
\int_{c}^{d}f\left( x\right) dx  \notag \\
&\leq &\frac{f_{-}^{\prime }\left( b\right) \left[ \left( b-d\right)
^{2}+\left( b-d\right) \left( b-c\right) +\left( b-c\right) ^{2}\right] }{%
6\left( b-a\right) }  \notag \\
&&-\frac{f_{+}^{\prime }\left( a\right) \left[ \left( d-a\right) ^{2}+\left(
d-a\right) \left( c-a\right) +\left( c-a\right) ^{2}\right] }{6\left(
b-a\right) }.  \notag
\end{eqnarray}
\end{theorem}

\begin{proof}
Since $f$ is convex, then for $a.e.$ $x\in \left[ a,b\right] $, we have (by (%
\ref{2.9})) that 
\begin{equation}
\left( \frac{a+b}{2}-x\right) f^{\prime }\left( x\right) \leq \frac{1}{b-a}%
\int_{a}^{b}f\left( t\right) dt-f\left( x\right) .  \label{4.2}
\end{equation}
Integrating (\ref{4.2}) on $\left[ c,d\right] $ we deduce 
\begin{equation}
\frac{1}{d-c}\int_{c}^{d}\left( \frac{a+b}{2}-x\right) f^{\prime }\left(
x\right) dx\leq \frac{1}{b-a}\int_{a}^{b}f\left( t\right) dt-\frac{1}{d-c}%
\int_{c}^{d}f\left( x\right) dx.  \label{4.3}
\end{equation}
Since 
\begin{eqnarray*}
&&\frac{1}{d-c}\int_{c}^{d}\left( \frac{a+b}{2}-x\right) f^{\prime }\left(
x\right) dx \\
&=&\frac{1}{d-c}\left[ \left( \frac{a+b}{2}-d\right) f\left( d\right)
-\left( \frac{a+b}{2}-c\right) f\left( c\right) +\int_{c}^{d}f\left(
x\right) dx\right] 
\end{eqnarray*}
then by (\ref{4.3}) we deduce the first part of (\ref{4.1}).

Using (\ref{2.11}), we may write for any $x\in \left[ a,b\right] $ that 
\begin{equation}
\frac{1}{b-a}\int_{a}^{b}f\left( t\right) dt-f\left( x\right) \leq \frac{1}{%
2\left( b-a\right) }\left[ \left( b-x\right) ^{2}f_{-}^{\prime }\left(
b\right) -\left( x-a\right) ^{2}f_{+}^{\prime }\left( a\right) \right] .
\label{4.4}
\end{equation}
Integrating (\ref{4.4}) on $\left[ c,d\right] ,$ we deduce 
\begin{eqnarray}
&&\frac{1}{b-a}\int_{a}^{b}f\left( t\right) dt-\frac{1}{d-c}%
\int_{c}^{d}f\left( x\right) dx  \label{4.5} \\
&\leq &\frac{1}{2\left( b-a\right) }\left[ f_{-}^{\prime }\left( b\right) 
\frac{1}{d-c}\int_{c}^{d}\left( b-x\right) ^{2}dx-f_{+}^{\prime }\left(
a\right) \frac{1}{d-c}\int_{c}^{d}\left( x-a\right) ^{2}dx\right] .  \notag
\end{eqnarray}
Since 
\begin{equation*}
\frac{1}{d-c}\int_{c}^{d}\left( b-x\right) ^{2}dx=\frac{\left( b-d\right)
^{2}+\left( b-d\right) \left( b-c\right) +\left( b-c\right) ^{2}}{3}
\end{equation*}
and 
\begin{equation*}
\frac{1}{d-c}\int_{c}^{d}\left( x-a\right) ^{2}dx=\frac{\left( d-a\right)
^{2}+\left( d-a\right) \left( c-a\right) +\left( c-a\right) ^{2}}{3},
\end{equation*}
then by (\ref{4.5}) we deduce the second part of (\ref{4.1}).
\end{proof}

\begin{remark}
If we choose $f\left( x\right) =x^{p},\;p\in \left( -\infty ,0\right) \cup
\lbrack 1,\infty )\backslash \left\{ -1\right\} $ or $f\left( x\right) =%
\frac{1}{x}$ \ or even $f\left( x\right) =-\ln x$, $x\in \left[ a,b\right]
\subset \left( 0,\infty \right) ,$ in the above inequalities, then a great
number of interesting results for $p-$logarithmic, logarithmic and identric
means may be obtained. We leave this as an exercise to the interested reader.
\end{remark}

\section{Applications for P.D.F.s}

Let $X$ be a random variable with the \textit{probability density function} $%
f:\left[ a,b\right] \subset \mathbb{R\rightarrow R}_{+}$ and with \textit{%
cumulative distribution function }$F\left( x\right) =\Pr \left( X\leq
x\right) .$

The following theorem holds.

\begin{theorem}
\label{t5.1}If $f:\left[ a,b\right] \subset \mathbb{R\rightarrow R}_{+}$ is
monotonically increasing on $\left[ a,b\right] $, then we have the
inequality: 
\begin{eqnarray}
&&\frac{1}{2}\left[ \left( b-x\right) ^{2}f_{+}\left( x\right) -\left(
x-a\right) ^{2}f_{-}\left( x\right) \right]  \label{5.1} \\
&\leq &b-E\left( X\right) -\left( b-a\right) F\left( x\right)  \notag \\
&\leq &\frac{1}{2}\left[ \left( b-x\right) ^{2}f_{-}\left( b\right) -\left(
x-a\right) ^{2}f_{+}\left( a\right) \right]  \notag
\end{eqnarray}
for any $x\in \left( a,b\right) ,$ where $f_{-}\left( \alpha \right) $ means
the left limit in $\alpha $ while $f_{+}\left( \alpha \right) $ means the
right limit in $\alpha $ and $E\left( X\right) $ is the expectation of $X.$%
\newline
The constant $\frac{1}{2}$ is sharp in both inequalities.\newline
The second inequality also holds for $x=a$ of $x=b.$
\end{theorem}

\begin{proof}
Follows by Theorem \ref{t1} and \ref{t2} applied for the convex cdf function 
$F\left( x\right) =\int_{a}^{x}f\left( t\right) dt,\;x\in \left[ a,b\right] $
and taking into account that 
\begin{equation*}
\int_{a}^{b}F\left( x\right) dx=b-E\left( X\right) .
\end{equation*}
\end{proof}

Finally, we may state the following corollary in estimating the probability $%
\Pr \left( X\leq \frac{a+b}{2}\right) .$

\begin{corollary}
\label{c5.2}With the above assumptions, we have 
\begin{eqnarray}
&&b-E\left( X\right) -\frac{1}{8}\left( b-a\right) ^{2}\left[ f_{-}\left(
b\right) -f_{+}\left( a\right) \right]  \label{4.2} \\
&\leq &\Pr \left( X\leq \frac{a+b}{2}\right)  \notag \\
&\leq &b-E\left( X\right) -\frac{1}{8}\left( b-a\right) ^{2}\left[
f_{+}\left( \frac{a+b}{2}\right) -f_{-}\left( \frac{a+b}{2}\right) \right] .
\notag
\end{eqnarray}
\end{corollary}

\section{Applications for $HH-$Divergence}

Assume that a set $\chi $ and the $\sigma -$finite measure $\mu $ are given.
Consider the set of all probability densities on $\mu $ to be 
\begin{equation}
\Omega :=\left\{ p|p:\Omega \rightarrow \mathbb{R},\;p\left( x\right) \geq
0,\;\int_{\chi }p\left( x\right) d\mu \left( x\right) =1\right\} .
\label{6.1}
\end{equation}
Csisz\'{a}r's $f-$divergence is defined as follows \cite{11b} 
\begin{equation}
D_{f}\left( p,q\right) :=\int_{\chi }p\left( x\right) f\left[ \frac{q\left(
x\right) }{p\left( x\right) }\right] d\mu \left( x\right) ,\;p,q\in \Omega ,
\label{6.2}
\end{equation}
where $f$ is convex on $\left( 0,\infty \right) $. It is assumed that $%
f\left( u\right) $ is zero and strictly convex at $u=1.$ By appropriately
defining this convex function, various divergences are derived.

In \cite{12b}, Shioya and Da-te introduced the generalised Lin-Wong $f-$%
divergence $D_{f}\left( p,\frac{1}{2}p+\frac{1}{2}q\right) $ and the
Hermite-Hadamard $\left( HH\right) $ divergence 
\begin{equation}
D_{HH}^{f}\left( p,q\right) :=\int_{\chi }\frac{p^{2}\left( x\right) }{%
q\left( x\right) -p\left( x\right) }\left( \int_{1}^{\frac{q\left( x\right) 
}{p\left( x\right) }}f\left( t\right) dt\right) d\mu \left( x\right)
,\;p,q\in \Omega ,  \label{6.3}
\end{equation}
and, by the use of the Hermite-Hadamard inequality for convex functions,
proved the following basic inequality 
\begin{equation}
D_{f}\left( p,\frac{1}{2}p+\frac{1}{2}q\right) \leq D_{HH}^{f}\left(
p,q\right) \leq \frac{1}{2}D_{f}\left( p,q\right) ,  \label{6.4}
\end{equation}
provided that $f$ is convex and normalised, i.e., $f\left( 1\right) =0.$

The following result in estimating the difference 
\begin{equation*}
D_{HH}^{f}\left( p,q\right) -D_{f}\left( p,\frac{1}{2}p+\frac{1}{2}q\right)
\end{equation*}
holds.

\begin{theorem}
\label{t6.1}Let $f:[0,\infty )\rightarrow \mathbb{R}$ be a convex function
and $p,q\in \Omega .$ Then we have the inequality: 
\begin{eqnarray}
0 &\leq &\frac{1}{8}\left[ D_{f_{+}^{\prime }\cdot \left| \frac{\cdot +1}{2}%
\right| }\left( p,q\right) -D_{f_{-}^{\prime }\cdot \left| \frac{\cdot +1}{2}%
\right| }\left( p,q\right) \right]  \label{6.5} \\
&\leq &D_{HH}^{f}\left( p,q\right) -D_{f}\left( p,\frac{1}{2}p+\frac{1}{2}%
q\right)  \notag \\
&\leq &\frac{1}{8}D_{f_{-}^{\prime }\cdot \left( \cdot -1\right) }\left(
p,q\right) .  \notag
\end{eqnarray}
\end{theorem}

\begin{proof}
Using the double inequality 
\begin{eqnarray*}
0 &\leq &\frac{1}{8}\left[ f_{+}^{\prime }\left( \frac{a+b}{2}\right)
-f_{-}^{\prime }\left( \frac{a+b}{2}\right) \right] \left| b-a\right|  \\
&\leq &\frac{1}{b-a}\int_{a}^{b}f\left( t\right) dt-f\left( \frac{a+b}{2}%
\right)  \\
&\leq &\frac{1}{8}\left[ f_{-}\left( b\right) -f_{+}^{\prime }\left(
a\right) \right] \left( b-a\right) 
\end{eqnarray*}
for the choices $a=1$, $b=\frac{q\left( x\right) }{p\left( x\right) },$ $%
x\in \chi ,$ multiplying with $p\left( x\right) \geq 0$ and integrating over 
$x$ on $\chi $ we get 
\begin{eqnarray*}
0 &\leq &\frac{1}{8}\int_{\chi }\left[ f_{+}^{\prime }\left( \frac{p\left(
x\right) +q\left( x\right) }{2p\left( x\right) }\right) -f_{-}^{\prime
}\left( \frac{p\left( x\right) +q\left( x\right) }{2p\left( x\right) }%
\right) \right] \left| q\left( x\right) -p\left( x\right) \right| d\mu
\left( x\right)  \\
&\leq &D_{HH}^{f}\left( p,q\right) -D_{f}\left( p,\frac{1}{2}p+\frac{1}{2}%
q\right)  \\
&\leq &\frac{1}{8}\int_{\chi }\left[ f_{-}^{\prime }\left( \frac{q\left(
x\right) }{p\left( x\right) }\right) -f_{+}^{\prime }\left( 1\right) \right]
\left( q\left( x\right) -p\left( x\right) \right) d\mu \left( x\right) ,
\end{eqnarray*}
which is clearly equivalent to (\ref{6.5}).
\end{proof}

\begin{corollary}
\label{c6.2}With the above assumptions and if $f$ is differentiable on $%
\left( 0,\infty \right) ,$ then 
\begin{equation}
0\leq D_{HH}^{f}\left( p,q\right) -D_{f}\left( p,\frac{1}{2}p+\frac{1}{2}%
q\right) \leq \frac{1}{8}D_{f^{\prime }\cdot \left( \cdot -1\right) }\left(
p,q\right) .  \label{6.6}
\end{equation}
\end{corollary}

\end{document}